\documentclass[11pt,leqno]{article}

\usepackage{amsmath,amssymb}
\usepackage{epic}
\usepackage{epsfig}

\begin{document}

\newtheorem{defi}{Definition}[section]
\newtheorem{rem}{Remark}[section]
\newtheorem{prop}{Proposition}[section]
\newtheorem{lem}{Lemma}[section]
\newtheorem{theo}{Theorem}[section]
\newtheorem{cor}{Corollary}[section]
\newtheorem{conc}{Conclusion}[section]
\newtheorem{claim}{Claim}[section]

\author{Pierre Germain}

\title{On the existence of smooth self-similar blow-up profiles for the wave-map equation}

\maketitle

\maketitle

\begin{abstract}
Consider the equivariant wave map equation from Minkowski space to a rotationnally symmetric manifold $N$ which has an equator (example: the sphere).

In dimension $3$, this article presents a necessary and sufficient condition on $N$ for the existence of a smooth self-similar blow up profile.

More generally, we study the relation between

- the minimizing properties of the equator map for the Dirichlet energy corresponding to the (elliptic) harmonic map problem

- and the existence of a smooth blow-up profile for the (hyperbolic) wave map problem.

This has several applications to questions of regularity and uniqueness for the wave map equation.
\end{abstract}

\section{Introduction}

\subsection{The geometric setting}

\label{lezard}

We will consider in the following equations on a map $u$, from the Euclidean space $\mathbb{R}^{d}$, with $d \geq3$ (which corresponds to the supercritical regime for harmonic maps and wave maps) to a manifold $N$ embedded in $\mathbb{R}^K$.

\subsubsection{Spherical symmetry and equator}

The first assumption that we will make is one of symmetry: on $N$, one can take coordinates (which may be periodic in $\phi$)
$$
(\phi,\chi) \in {\mathbb{R}} \times \mathbb{S}^{M-1}
$$
($M$ an integer) in which the metric reads
$$
d\phi^2 + g(\phi)^2 d\chi^2\,\,,
$$
for a smooth odd function $g$ such that $g'(0)=1$. 

\bigskip

We now impose that $N$ has an equator; by this we simply mean $\phi^*>0$ such that
$$
g'(\phi^*) = 0 \,\,,
$$
and minimal with this property.

\subsubsection{Equivariance}

We further restrict $u$ to the class of equivariant maps, that is, $(r,\omega)$ denoting polar coordinates on $\mathbb{R}^d$,
\begin{equation*}
\begin{split}
u = R_\chi \phi\;\overset{def}{:}\; & {\mathbb{R}}^{d} \longrightarrow N \\
& (r,\omega) \longrightarrow (\phi(r)\,,\chi(\omega)) \,\,.
\end{split}
\end{equation*}
Here, $\chi : \mathbb{S}^{d-1} \rightarrow \mathbb{S}^{M-1}$ ($M$ is an integer) is an eigenmap, that is a map satisfying, for a constant $k$
\begin{equation*}
|\nabla \chi|^2 = k \;\;\;\;\;\mbox{and}\;\;\;\;\; \Delta_{\mathbb{S}^{d-1}} \chi + k \chi = 0\,\,.
\end{equation*}
This can be achieved for some large enough $M$, provided
$$
k = \ell (\ell + d -2) \;\;\;\;\mbox{with $\ell \in \mathbb{N}$}
$$
see~\cite{CSTZ}.

\bigskip

Of course, the most simple example of an eigenmap is provided by $\chi(\omega) = \omega$, mapping $\mathbb{S}^{d-1}$ to itself, in which case $k = d-1$.

\bigskip

We call equator map the map corresponding to
$$
\phi \equiv \phi^* \,\,;
$$
it actually depends on $\chi$, that is the equivariance class which is considered.

\subsection{Harmonic maps and the equator map}

\subsubsection{The boundary value problem}

Harmonic maps from the ball $B = B_{{\mathbb{R}}^d}(0,1)$ of the Euclidean space $\mathbb{R}^{d}$, to a manifold $N$ embedded in $\mathbb{R}^K$, with Dirichlet boundary data $h$, satisfy the equation
$$
\left\{ \begin{array}{ll}
 \Delta u(x) \perp T_{u(x)} N & \mbox{on $B$} \\
 u = h & \mbox{on $\partial B$} \end{array} \right.
$$
(here, $T_{u(x)} N$ is the tangent space, at $u(x)$, of $N$, and orthogonality is to be understood in the sense of $\mathbb{R}^{K}$).
In the equivariant setting discussed above, this Dirichlet problem becomes
\begin{equation}
\label{dirichletequivariant}
(D_\alpha) \;\;\;\;
\left\{ \begin{array}{ll}
 -\phi_{rr} - \frac{d-1}{r} \phi_r + \frac{k}{r^2} g(\phi) g'(\phi) = 0 & \mbox{on $[0,1]$} \\
 \phi(1) = \alpha \,\,, & \end{array} \right.
\end{equation}
where $\alpha$ is a given real number.

\subsubsection{The energy and the equator map}

Harmonic maps are critical points of the energy, which reads in the equivariant setting
$$
E_e(\phi) = \int_0^1 \left( \phi_r^2 + \frac{k}{r^2} g^2(\phi) \right) r^{d-1} dr\,\,.
$$
Solutions to the Dirichlet problem $(D_\alpha)$ are critical points of $E_e$ in the class
$$
F_\alpha = \left\{ \phi:[0,1]\rightarrow {\mathbb{R}}\;,\; \phi_r \in L^2(r^{d-1}dr) \;,\;\phi(1) = \alpha \right\}\,\,.
$$
It can easily be seen that the equator map $\phi \equiv \phi^*$ is a solution of the Dirichlet problem with $\alpha = \phi^*$. A crucial point in this article will be the minimizing properties of $\phi^*$. 
\begin{itemize}
\item
For a full answer to the question whether $\phi^*$ minimizes locally $E_e$ over $F_{\phi^*}$, see Proposition~\ref{local}.
\item
Whether $\phi^*$ minimizes globally or not $E_e$ over $F_{\phi^*}$ is a more difficult question. The answer is known in the case where $\chi = Id$ (covariant case) and $N$ is an ellipse embedded in ${\mathbb{R}}^{d+1}$ given by ($a \in [0,1]$)
$$
y_1^2 + \dots + y_d^2 + \frac{y_{d+1}^2}{a^2} = 1 \,\,,
$$
and it is due to J\"ager and Kaul~\cite{JK}, Baldes~\cite{Ba}, H\'elein~\cite{He}. 
These authors prove the following criterion: $\phi^*$ is the unique global minimizer of $E_e$ if and only if
$$
a^2 \geq \frac{4(d-1)}{(d-2)^2}\,\,.
$$
\end{itemize}

\subsection{Wave maps}

This paper addresses only some of the questions related to the wave map Cauchy problem. Other aspects, like the case of dimension $2$, or optimal well-posedness results for wave maps without symmetries, are discussed in the review by Krieger~\cite{K}.

\subsubsection{The initial value problem}

The Cauchy problem for wave maps from Minkowski space $\mathbb{R}^{d+1}$ to a manifold $N$ embedded in $\mathbb{R}^K$
reads
$$
(WM)\;\;\;\;
\left\{ \begin{array}{l}
 \partial_t^2 u(x,t) - \Delta u(x,t) \perp T_{u(x,t)} N \\
 u(t=0) = u_0 \\
 u_t(t=0) = u_1 \,\,. \end{array} \right.
$$
The weak formulation reads, for $(u\,,\,u_t) \in L^\infty_t (\dot{H}^1 \times L^2)$, $u(x,t) \in N$ for almost all $(x\,,\,t)$,
\begin{equation*}
\begin{split}
&\int_{\mathbb{R}^{d+1}} \left[ \langle u_t\,,\, f_t \rangle - \langle \nabla_x u \,,\, \nabla_x f \rangle \right] \,dx\,dt = 0 \\
& \;\;\;\;\;\;\;\;\;\;\;\;\;\;\;\;\;\;\;\;\;\;\;\;\mbox{for $f$ such that $f_t \,, \,\nabla_x f \in \L^2$ and for almost all $(x\,,\,t)$, $f(x,t) \in T_{u(x,t)} N$} \\
&(u \,,\,u_t) \overset{t \rightarrow 0}{\longrightarrow} (u_0\,,\,u_1) \;\;\mbox{weakly in $\dot{H}^1 \times L^2$}\,\,.
\end{split}
\end{equation*}
We can rewrite these equations in the equivariant setting. We obtain the following equation, which is to be understood in $\mathcal{S}'$ if one considers weak solutions.
$$
(EWM)\;\;\;\;
\left\{ \begin{array}{l}
 \phi_{tt} -\phi_{rr} - \frac{d-1}{r} \phi_r + \frac{k}{r^2} g(\phi) g'(\phi) = 0 \\
 \phi(t=0) = \phi_0 \\
 \phi_t(t=0) = \phi_1 \,\,. \end{array} \right.
$$

\subsubsection{Formation of singularities: the blow up profiles}

It is a basic problem to understand whether the wave map equation can develop singularities if the initial data are smooth. The approach followed by Shatah~\cite{Sh}, Shatah and Tahvildar-Zadeh~\cite{STZ}, and Cazenave, Shatah and Tahvildar-Zadeh~\cite{CSTZ} consists in building up smooth blow up profiles

\begin{defi}[blow up profiles]
A (self-similar) blow up profile is a function $\psi$ such that
$$
\phi(r,t) = \psi \left( \frac{r}{t} \right)
$$
solves the wave-map equation.
\end{defi}

 The profile $\psi(\rho)$ has to be defined for $\rho \leq 1$, and must have a smooth behavior at 1. Then one can extend the data $(\psi, r \psi')$, defined on $B_{\mathbb{R}^d}(0,1)$, smoothly to 
$(\widetilde{\phi}_0\,,\,\widetilde{\phi}_1)$, defined on $\mathbb{R}^d$. The solution to $(EWM)$ with data $(\widetilde{\phi}_0\,,\,\widetilde{\phi}_1)$ equals
$$
\psi \left( \frac{r}{1-t} \right)
$$
inside the backward light cone, with basis $B$ and vertex $(x=0,t=1)$. In particular, a singularity appears at $(x=0,t=1)$, even though the data were smooth.

\bigskip

The PDE satisfied by $\phi$ turns into an ODE for $\psi$. More precisely,

\begin{prop}
\label{PGProp}
Suppose that $\psi$ is an integrable function, smooth at $0$. Then
$R_\chi \phi$ solves $(WM)$ for $t>0$ if and only if $\psi$ satisfies the equation
\begin{equation}
\label{eqpsi}
\rho^2 (1-\rho^2) \psi''(\rho) + \rho ((d-1)- 2 \rho^2) \psi'(\rho) - k g(\psi) g'(\psi) = 0
\end{equation}
in the classical sense on $(0,1) \cup (1,\infty)$, and furthermore $\psi(1-) = \psi(1+)$.
\end{prop}

This proposition is proved in Section~\ref{sectionODE}. 

\bigskip

A first approach to solving for $\psi$ is to handle the above singular ODE. Another approach is variational, it was developed by the authors mentioned above. These authors observe that the equation satisfied by $\psi$ is the harmonic map equation from the hyperbolic ball to $N$. Thus it can be solved by minimizing
$$
E_h(\psi) = \int_{B} \left( \psi_r^2 + \frac{k}{r^2(1-r^2)} \left[ g^2(\psi) - g^2(\phi^*) \right] \right) \frac{r^{d-1} dr}{(1-r^2)^{(d-3)/2}}\,\,,
$$
on
$$
G = \left\{ \psi : [0,1] \rightarrow {\mathbb{R}}\;,\; \psi_r \in L^2\left(\frac{r^{d-1}\,dr}{(1-r^2)^{\frac{d-3}{2}}}\right)\;,\;\psi(1) = \phi^* \right\}\,\,.
$$
The functional $E_h$ corresponds to the Dirichlet energy for an equivariant map from the hyperbolic ball to $N$, except $g^2(\psi)$ has been replaced by $\left[ g^2(\psi) - g^2(\phi^*) \right]$. This renormalization imposes 
the boundary condition $\psi = \phi^*$ for $\rho = 1$ - we will come back to this in the following.

The following theorem summarizes results proved in~\cite{CSTZ}.

\begin{theo} 
\label{ehnegative}
If 
$$
\inf_{\psi \in G} E_h(\psi) < 0\,\,,
$$
there exists a smooth $\psi : [0,1] \rightarrow [0,R^*]$ such that 
\begin{itemize}
\item $\psi(0) = 0$.
\item $\psi(1) = \phi^*$.
 \item  $\phi = \psi \left( \frac{r}{t} \right)$ solves $(EWM)$ (where it is defined).
\item $u = R_\chi \phi$ is smooth where it is defined.
\end{itemize}
\end{theo}

It should also be mentioned that self-similar blow-up profiles as given by the above theorem provide examples of non uniqueness for $(EWM)$: see~\cite{STZ}~\cite{CSTZ}. We will use a similar argument in Theorem~\ref{dim3} for data equal to the equator map.

In the absence of such a self-similar blow-up profile, there does not seem to be known examples of non-uniqueness.

Finally, multiple solutions for the above variational problem are constructed in Jungen~\cite{J}.

\subsubsection{Energy conservation and weak solutions}

A smooth solution of $(EWM)$ satisfies the energy equality (we only give here the version which is centered at $0$ in space)
\begin{equation}
\label{energyeq}
\mbox{if $R>T$,}\;\;\;\;E(T, R - T,\phi) + \operatorname{flux}(0,T, R ,\phi) = E(0, R ,\phi)
\end{equation}
with the following notations
\begin{equation*}
\begin{split}
& E(t,R,\phi) \overset{def}{=} \int_0^R \left( \phi_r^2(t) + \frac{k}{r^2} g(\phi(t))^2 + \phi_t^2(t) \right) r^{d-1} \,dr \\
& \operatorname{flux}(0,T,R,\phi) \overset{def}{=} \int_{S_{0,T,R}} \left( \frac{1}{2} \phi_t^2 + \frac{k}{2r^2} g(\phi)^2 + \frac{1}{2} |\nabla \phi|^2 - \phi_t \nabla \phi \cdot n \right) d\sigma \,\,,
\end{split}
\end{equation*}
where $S_{0,T,R}$ is the surface
$$
S_{0,T,R} = \{ (r,s)\;\;\mbox{such that $0 \leq t \leq T$ and $r = R-s$} \}
$$
(endowed with the metric coming from the natural embedding $S \subset \mathbb{R}^{d+1}$), 
and where $n$ is the outside normal of $B(0,R - s)$. Notice that $\operatorname{flux}(0,T,R,\phi)$ is a positive quantity.

\bigskip

In the case where the data are of finite energy, one can use energy conservation and a compactness argument to build up global weak solutions.

\begin{theo}[Shatah~\cite{Sh}, Freire~\cite{Fr}]
Assume that the target manifold $N$ is homogeneous and compact, and 
take data such that $\nabla u_0 \in L^2_{loc}$ and $u_1 \in L^2_{loc}$, i.e. for any $R>0$,
$$
\int_0^R \left( (\phi_0)_r^2 + \frac{k}{r^2} g(\phi_0)^2 \right) r^{d-1}\,dr + \int_0^R \phi_1^2 r^{d-1}\,dr < \infty \,\,.
$$
Then there exists a global solution of $(WM)$
$$
(u,u_t) \in L^\infty_t(\mathbb{R},\dot{H}^1_{loc}) \times L^\infty_t(\mathbb{R},L^2_{loc}) \,\,.
$$
\end{theo}

If one transposes this theorem in the equivariant setting, the following energy inequality can be proved to hold
\begin{equation}
\label{energyineq}
\mbox{if $R>T$,}\;\;\;\;E(T, R - T,\phi) + \operatorname{flux}(0,T, R ,\phi) \leq E(0, R ,\phi)\,\,.
\end{equation}
Indeed, it suffices to follow the steps of Shatah and Freire, who use a penalisation method. Since we consider equivariant functions, passing to the limit in the penalised problems is easy except for $x=0$; this is due to the energy bound which gives strong control anywhere else. Thus one can pass to the limit in~(\ref{energyeq}) and obtain~(\ref{energyineq}).

\section{Statement of the results}

\subsection{Local aspects}

\begin{prop}
\label{local}
The three following facts are equivalent
\medskip

(i) $k g(\phi^*) g''(\phi^*) \geq -\frac{(d-2)^2}{4}$
\medskip

(ii) The second variation $\displaystyle \delta^2 E_e(\phi^*)(w,w) \geq 0$ for $w \in F_{\phi^*}$.
\medskip

(iii) The second variation $\displaystyle \delta^2 E_h(\phi^*)(w,w) \geq 0$ for $w \in G$.
\bigskip

Besides, if we exclude the limit case where $k g(\phi^*) g''(\phi^*) = -\frac{(d-2)^2}{4}$, $(i)$ becomes equivalent to

\medskip

(iv) The wave map equation with initial data equal to the equator map is linearly stable (in a sense to be made precise in Section~\ref{linearlystable}).
\end{prop}

The proof of this proposition essentially consists in putting together already known results. It is interesting to notice that they all rely on different versions of Hardy's inequality.

This proof will be given in Section~\ref{prooflocal}.

\subsection{The case of dimension $3$}

In dimension $3$, we are able to analyze globally the relation between minimality of $\phi^*$ for $E_e$ and properties of the wave-map equation.

However, we need to make some more hypotheses on the geometry of $N$. We do not claim that they are optimal, but they are simple, of sufficient generality, and will enable us to state our results in an elegant way. Recall that $g$ is an odd smooth function. The further assumptions in dimension 3 are that
\begin{itemize}

\item {\bf either} the $\phi$ coordinate is not periodic; $\phi^*$ and $- \phi^*$ are the only zeros of $g'$; $g$ is positive and decreasing on $(\phi^*,\infty)$.

\item {\bf or} the $\phi$ coordinate is $4 \phi^*$ periodic, and $g (\phi^* + \cdot)$ is even. (Thus the manifold $N$ has the symmetries of an ellipsoid with all semi-principal axes but one of the same length)

\end{itemize}

Either of these assumptions ensures that
\begin{itemize}
\item Up to the symmetries of $N$, $\phi^*$ is unique with the property that $g'(\phi^*) = 0$.
\item Up to the symmetries of $N$, $0$ is unique with the property that $g'(\phi^*)=0$.
\item Up to the symmetries of $N$, 
\begin{equation}
\label{etoile}
g(\phi)^2 < g(\phi^*)^2 \;\;\;\mbox{for any}\;\; \phi \neq \phi^* \,\,.
\end{equation}
\end{itemize}

We can now state the theorem.

\begin{theo}
\label{dim3}
If $d=3$, and the above assumptions are satisfied, the three following facts are equivalent
\medskip

(i) The equator map is the unique global minimizer of $E_e$ on $F_{\phi^*}$.
\medskip

(ii) The only weak solution to $(EWM)$ satisfying the energy inequality~(\ref{energyineq}), and with initial data
$$\phi_0 = \phi^* \;\;\;\;\;\;\phi_1 = 0$$
is identically equal to the equator map $\phi^*$.
\medskip

(iii) There does not exist a non-zero blow-up profile $\psi$ such that $R_\chi \psi \in \mathcal{C}^\infty (\bar B_{\mathbb{R}^d}(0,1))$.

\end{theo}

This theorem will be proved in Section~\ref{prodim3}. Let us mention quickly the main ideas of the proof.
\begin{itemize}
\item Starting from a smooth blow up profile $\psi$, it is possible to construct, using the equation it satisfies, a solution corresponding to data $(\phi_0 \,,\, \phi_1) = (\phi^*\,,\,0)$ and non identically equal to $\phi^*$.
\item If the equator map minimizes the Dirichlet energy, it is intuitively clear that the equator map is the unique solution corresponding to data $(\phi^* \,,\, 0)$ and satisfying the energy inequality. Indeed, any other solution would somehow have to make the energy grow, which is not allowed if the energy inequality is satisfied.
\item Finally, if the equator map does not minimize the Dirichlet energy, we use the variational approach: Theorem~\ref{ehnegative} gives a smooth solution.
\end{itemize}

Actually, the proof of Theorem~\ref{dim3} gives a little more than the statement above. In particular, we get the following very strong instability result, for data equal to the equator map, in case $(i)$ above does not hold.

\begin{prop}
Suppose that the equator map $\phi^*$ is not the unique global minimizer of $E_e$. Then there exists a smooth profile $\psi: [0,1] \rightarrow \mathbb{R}$ such that for any $T \geq 0$,
$$
\phi (r,t) = \left\{ \begin{array}{ll}
\psi\left( \frac{r}{t-T} \right) & \mbox{if $t \geq T$ and $r \leq t-T$} \\
\phi^* & otherwise \end{array}
\right.
$$
solves $(EWM)$ with data $(\phi_0 \,,\, \phi_1) = (\phi^*\,,\,0)$, and satisfies the energy equality~(\ref{energyeq}).
\end{prop}

Finally, it is natural to ask whether non-smooth profiles may exist even if $(i)$ in Theorem~\ref{dim3} holds. In other words: can one replace in the statement of $(iii)$ ``smooth'' by ``$H^s$'' for some $s$? 
(Notice that the belonging of $R_\chi \psi$ to any functional space can be easily defined using the embedding of $N$ in $\mathbb{R}^K$.)

\begin{prop}
(i) Suppose that $\psi$ solves~(\ref{eqpsi}) on $(0,1)$, and that
$$
R_\chi \psi \in \dot{H}^{3/2}(B_{\mathbb{R}^3}(0,1))\,\,.
$$
Then $R_\chi \psi \in \mathcal{C}^\infty(\bar B_{\mathbb{R}^3}(0,1))$.

\bigskip

(ii) With only the assumptions made in Section~\ref{lezard} on $g$, there always exists a non zero $\psi$ such that
$R_\chi \psi$ smooth near $0$,
$$
R_\chi \psi \in \dot{B}^{3/2}_{2,\infty}(B_{\mathbb{R}^3}(0,1))\,\,,
$$
(see~\cite{Ge} for a definition of the Besov space $\dot{B}^{3/2}_{2,\infty}$) and
$$
u(x,t) = R_\chi \psi \left( \frac{|x|}{t} \right)
$$
solves $(WM)$.

\end{prop}

This proposition is proved in Section~\ref{proofreg}.

Thus, for $\psi$, $\dot{H}^{3/2}$ and $\mathcal{C}^\infty$ regularity are equivalent. But if one goes down to $\dot{B}^{3/2}_{2,\infty}$ the situation becomes drastically different: there always exist (no matter whether $(i)$ in Theorem~\ref{dim3} holds true or not) a profile with this regularity.

\bigskip

We would like to shed some light on the distinction between $\dot{H}^{3/2}$ and $\dot{B}^{3/2}_{2,\infty}$ appearing above.

Notice first that the spaces $\dot{B}^{3/2}_{2,\infty} \times \dot{B}^{1/2}_{2,\infty}$ or $\dot{H}^{3/2} \times \dot{H}^{1/2}$ are at the scaling of the equation for the initial data.
Global well-posedness for small data in $\dot{H}^{3/2}\times \dot{H}^{1/2}$ for $(WM)$ is due to Tao~\cite{Ta}. Global well posedness in the equivariant setting for a space similar to $\dot{B}^{3/2}_{2,\infty}\times \dot{B}^{1/2}_{2,\infty}$ has been proved by the author~\cite{Ge}.

Now let us point out a fundamental difference between considering solutions in the Sobolev space or in the Besov space. Space-time scaling invariant spaces corresponding to the spaces for the data mentioned above are
$$
(u,u_t) \in L^\infty_t \left( \dot{B}^{3/2}_{2,\infty} \times \dot{B}^{1/2}_{2,\infty}\right)\;\;\;\;\mbox{or}\;\;\;\; L^\infty_t \left(\dot{H}^{3/2} \times \dot{H}^{1/2}\right)\,\,.
$$

Consider a smooth profile $\psi$ and
$$
u(x,t) = R_\chi \psi \left( \frac{|x|}{t} \right) \,\,.
$$

Then (at least locally in space)
\begin{equation}
\label{besov}
(u,u_t) \in L^\infty \left([-1,1] \,,\,\dot{B}^{3/2}_{2,\infty} \times \dot{B}^{1/2}_{2,\infty}\right)
\end{equation}
but
\begin{equation}
\label{sobolev}
(u,u_t) \notin L^\infty \left([-1,1] \,,\,\dot{H}^{3/2} \times \dot{H}^{1/2} \right)\,\,.
\end{equation}
Thus in some sense one ``sees'' the blow up in the Sobolev setting~(\ref{sobolev}) but not in the Besov space setting~(\ref{besov}).

\subsection{The case of dimension $d \geq 4$}

In the case of dimension $d \geq 4$, we can only prove a weak version of Theorem~\ref{dim3}. 

\begin{theo}
\label{dimg4}
Consider the following assertions
\medskip

(i) The equator map is the unique global minimizer of $E_e$ on $F_{\phi^*}$.
\medskip

(ii) The only weak solution to $(EWM)$ satisfying the energy inequality~(\ref{energyineq}), and with initial data
$$\phi_0 = \phi^* \;\;\;\;\;\;\phi_1 = 0$$
is identically equal to the equator map $\phi^*$.
\medskip

(iii) There does not exist a blow-up profile $\psi$ such that $\psi(1) = \phi^*$ and $R_\chi \psi \in \mathcal{C}^\infty (\bar B_{\mathbb{R}^d}(0,1))$.
\medskip

There holds
$$
(i) \Rightarrow (ii) \Rightarrow (iii)\,\,.
$$
\end{theo}

The proof of this theorem is almost identical to the corresponding steps in Theorem~\ref{dim3}; for this reason we do not include it here.

The question of the regularity of the profile is much more involved in dimensions $d \geq 4$: see~\cite{CSTZ}; we plan to come back to this and other aspects of the problem for $d \geq 4$ in a subsequent paper. 

\section{Proof of Proposition~\ref{local}}

\label{prooflocal}

As mentioned above, the whole proof of Proposition~\ref{local} somehow relies on different versions of Hardy's inequality. It is interesting to recall here the whole space version, with the optimal constant
$$
\int_{\mathbb{R}^d} \frac{w^2}{x^2} \,dx \leq \frac{4}{(d-2)^2} \int_{\mathbb{R}^d} |\nabla w|^2 \,dx\,\,.
$$

\subsection{The second variation of $E_e$: proof of $(i) \Leftrightarrow (ii)$}

A small computation gives
$$
E_e (\phi^* + \epsilon w) = E_e(\phi^*) + \epsilon^2 \int_0^1 \left( w_r^2 + \frac{k}{r^2} g(\phi^*) g''(\phi^*) w^2 \right) r^{d-1} \,dr + O(\epsilon^3)\,\,.
$$
The conclusion follows due to the inequality (Baldes~\cite{Ba})
\begin{equation}
\label{HI1}
\mbox{if $\phi^* + \epsilon w \in F_{\phi^*}$,}\;\;\;\int_0^1 w^2 r^{d-3} dr \leq \frac{4}{(d-2)^2} \int_0^1 w_r^2 r^{d-1} dr\,\,,
\end{equation}
where the constant is optimal.

\subsection{The second variation of $E_h$: proof of $(i) \Leftrightarrow (iii)$}

By a small computation one gets (since $E_h(\phi^*) = 0$)
$$
E_h (\phi^* + \epsilon w) =  \epsilon^2 \int_0^1 \left( w_r^2 + \frac{k}{r^2(1-r^2)} g(\phi^*) g''(\phi^*) w^2 \right)
\frac{r^{d-1} dr}{(1-r^2)^{(d-3)/2}} + O(\epsilon^3) \,\,.
$$
The conclusion follows due to the inequality (Cazenave Shatah Tahvildar-Zadeh~\cite{CSTZ})
\begin{equation}
\label{HI2}
\mbox{if $\phi^* + \epsilon w \in G$,}\;\;\;
\int_0^1 w^2 \frac{r^{d-3}}{(1-r^2)^{\frac{d-1}{2}}} dr \leq \frac{4}{(d-2)^2} \int_0^1 w_r^2 \frac{r^{d-1}}{(1-r^2)^{\frac{d-3}{2}}} dr
\end{equation}
where the constant is optimal.

\subsection{Linear stability of the equator map: proof of $(i) \Leftrightarrow (iv)$}

\label{linearlystable}

The equation $(EWM)$ can be rewritten using the new unknown function $w = \phi - \phi^*$
\begin{equation}
\label{EWMw}
w_{tt} - w_{rr} - \frac{d-1}{r} w_r + \frac{k}{r^2} g(\phi^* + w) g'(\phi^* + w) = 0 \,\,.
\end{equation}
Let us now linearize this equation in the limit $w$ small:
\begin{equation}
\label{linearizedeq}
w_{tt} - w_{rr} - \frac{d-1}{r} w_r + \frac{k}{r^2} g(\phi^*) g''(\phi^*) w = 0 \,\,.
\end{equation}
This is a linear wave equation with a singular (inverse square) potential. Dispersion for this class of potentials has been studied by Burq, Planchon, Stalker and Tahvildar Zadeh~\cite{BPSTZ}, who prove the following result.

\begin{theo}[Burq, Planchon, Stalker, Tahvildar-Zadeh]
\label{theoBPSTZ}
The solution of the equation
$$
\left\{ \begin{array} {l}
w_{tt} - w_{rr} - \frac{d-1}{r} w_r + \frac{a}{r^2} w = 0 \\
w(t=0) = w_0 \\
w_t(t=0) = w_1 \\
\end{array} \right.
$$
satisfies the classical Strichartz estimates
\begin{equation}
\label{strichartz}
\|(-\Delta)^{\frac{\sigma}{2}} w \|_{L^p_t L^q_x} \leq C \left( \|w_0\|_{\dot{H}^{1/2}} + \|w_1\|_{\dot{H}^{-1/2}} \right)
\;\;\;\;\mbox{with $\frac{1}{p} + \frac{d-1}{2q} \leq \frac{d-1}{4}$ and $\sigma = \frac{1}{p} + \frac{d}{q} - \frac{d-1}{2}$}
\end{equation}
if
\begin{equation}
\label{criterion}
a > -\frac{(d-2)^2}{4} \,\,.
\end{equation}
\end{theo}

\bigskip

In the case where $a < - \frac{(d-2)^2}{4}$, the point spectrum of the operator $-\Delta + \frac{a}{x^2}$ extends to $-\infty$, see Remark 1.1 in Planchon, Stalker and Tahvildar-Zadeh~\cite{PSTZ}.

In particular, in this case, there are no Strichartz estimates.

\bigskip

What we mean by ``linearly stable'' in the statement of Proposition~\ref{local} is now clear: it means satisfying the Strichartz estimates given in Theorem~\ref{theoBPSTZ}.

Applying the criterion~(\ref{criterion}) to the linearized equation~(\ref{linearizedeq}), we get as expected the condition $(i)$ in Proposition~\ref{local}.

\begin{rem}[Nonlinear stability]
Of course a much more interesting question than linear stability is non-linear stability, that is, well posedness of~(\ref{EWMw}).

For a certain range of $k g(\phi^*) g''(\phi^*)$, a possible approach is the following: combine the Strichartz estimates with derivatives (that is, involving other regularities than what appears in~(\ref{strichartz})) of~\cite{BPSTZ} and the approach followed in~\cite{STZ} in order to get well posedness of~(\ref{EWMw}).

This approach, however, does not cover all the cases where linear stability holds. The difficulty is that, if $a$ is too close to $-\frac{(d-2)^2}{4}$, Strichartz estimates with derivatives are not available any more, and have to be replaced by estimates on
$$
\left\| (-\Delta)^{\sigma/2} \left( - \Delta + \frac{a}{x^2} \right)^{\sigma'} w \right\|_{Lp_t L^q_x}\,\,,
$$
that is the relevant spaces involve fractional powers of the operator $\displaystyle - \Delta + \frac{a}{x^2}$. But as opposed to standard Sobolev spaces, product laws between such spaces are not known and need to be established in order to be able to handle nonlinear settings.
\end{rem}

\begin{rem} 
Hardy's inequality did not appear explicitly in the present subsection, but it is one of the fundamental ingredients in the proof of Theorem~\ref{theoBPSTZ}.
\end{rem}

\section{Proof of Theorem~\ref{dim3}}

\label{prodim3}

\subsection{Equivalence of $(i)$ and $(ii)$}

\subsubsection{$(i) \Rightarrow (ii)$}

Let us suppose first that $(i)$ holds, and consider $\phi$ a weak solution of $(EWM)$ with data equal to the equator map, and satisfying the energy inequality~(\ref{energyineq}).

The equator map is smooth anywhere but in $0$. By finite speed of propagation and weak-strong uniqueness for $(EWM)$ (Struwe~\cite{St}), $u$ must agree with the equator map outside of the forward light cone with vertex at $(t,x) = (0,0)$
$$
C = \left\{ (t,x) \,,\,|x| \leq t \right\} \,\,.
$$ 
Let us write the energy equality~(\ref{energyeq}) for the equator map
\begin{equation}
\label{abc1}
E(T,R - T,\phi^*) + \operatorname{flux}(0,T,R,\phi^*) = E(0,R,\phi^*) \,\,,
\end{equation}
where we choose $2 T < R$ so that the flux term is computed on a surface which lies completely outside of the forward light cone $C$.

For the weak solution $\phi$, the energy inequality~(\ref{energyineq}) holds, which reads
\begin{equation}
\label{abc2}
E(T,R - T,\phi) + \operatorname{flux}(0,T,R,\phi) \leq E(0,R,\phi) \,\,.
\end{equation}

Comparing~(\ref{abc1}) and~(\ref{abc2}), and keeping in mind that $\phi$ and $\phi^*$ agree at time $0$ and outside of the light cone, one gets
$$
E(T,R - T,\phi^*) \geq E(T,R - T,\phi) \,\,,
$$
which implies
$$
\int_0^{R-T} \frac{k}{r^2} g(\phi^*)^2 \,r^2 dr \geq \int_0^{R-T} \left( \phi_r^2 (T) + \frac{k}{r^2} g(\phi(T))^2 \right)\,r^2 dr \,\,.
$$
Thus, since $\phi(T,R-T) = \phi^*$, the minimizing property of the equator map gives $\phi(T) = \phi^*$. In other words $(ii)$ holds true.

\subsubsection{$\overline{(i)} \Rightarrow \overline{ (ii)}$}

\label{cameleon}

Let us suppose that $(i)$ does not hold, and prove that $(EWM)$ with initial data equal to the equator map admits a solution which is not identically equal to the equator map.

The alternative solution that we build up will be a self-similar solution
$$
\phi (r,t) = \psi \left( \frac{r}{t} \right)
$$

Suppose that $\phi^*$ is not the unique global minimizer of $E_e$ in $F_{\phi^*}$. Then there exists $\widetilde{\phi} \in F_{\phi^*}$ such that 
$$
E_e(\widetilde{\phi}) \leq E_e(\phi^*) \,\,.
$$
Since the space dimension is $3$, $F_{\phi^*} = G$ and therefore it makes sense to compute
$$
E_h (\widetilde{\phi}) = E_e (\widetilde{\phi}) - E_e (\phi^*) + \int_0^1 \frac{k}{1-r^2}\left(g^2(\widetilde{\phi}) - g^2(\phi^*) \right)r^2\,dr \,\,.
$$
Since $\widetilde{\phi}$ is not identically $\phi^*$, and by condition~(\ref{etoile}),
$$
E_h (\widetilde{\phi}) < 0 \,\,.
$$
Hence Theorem~\ref{ehnegative} gives a self-similar profile $\psi$. Define the prolongation
$$
\bar \psi (\rho) = \left\{ \begin{array}{ll} \psi(\rho) & \mbox{if $\rho \leq 1$} \\ \phi^* & \mbox{if $\rho \geq 1$} \,\,.\end{array}
\right.
$$
and
$$
\phi(r,t) = \psi\left(\frac{r}{t} \right)\,\,.
$$
By proposition~\ref{PGProp}, $R_\chi \phi$ solves $(WM)$. It is easy to see that it takes the data
$$
\phi(t=0) = \phi^* \;\;\;\;\mbox{and}\;\;\;\; \phi_t (t=0) = 0 \,\,.
$$
To conclude that $(ii)$ is not true, it suffices to see that $\bar \psi$ satisfies the energy equality. This is precisely the statement of the next proposition, whose proof we postpone till Section~\ref{pec}.

\begin{prop}
\label{energyconserve}
$\bar \psi$, as defined above, satisfies the energy equality~(\ref{energyeq}).
\end{prop}

\subsection{Equivalence of $(i)$ and $(iii)$}

\subsubsection{$\overline{(i)} \Rightarrow \overline{(iii)}$}

Let us suppose that $(i)$ does not hold. Proceeding as in Section~\ref{cameleon}, we get a profile $\psi$ which is smooth on $[0,1]$. Thus $(iii)$ is not true.

\subsubsection{$\overline{(iii)} \Rightarrow \overline{(i)}$}

\label{subsubiii}

Let us suppose that $(iii)$ does not hold true, that is there exists a non zero blow up profile $\psi$ such that $R_\chi \psi \in \mathcal{C}^\infty( \bar B_{\mathbb{R}^2} (0,1) )$. Recall $\psi$ satisfies the ODE
$$
\rho^2 (1-\rho^2) \psi''(\rho) + 2 \rho (1-\rho^2) \psi'(\rho) - k g(\psi) g'(\psi) = 0
$$
Setting $\rho = 1$, we get $g(\psi(1)) g'(\psi(1)) = 0$. Thus, modulo the symmetries of the equation, either $\psi(1) = 0$ or $\psi(1) = \phi^*$. In order to show that the former case cannot occur, set
$$
M(\rho) = \frac{1}{2}(\rho^2 - \rho^4)\psi'(\rho)^2 - k g(\psi(\rho))^2
$$
and compute
$$
M'(\rho) = - \rho \psi'(\rho)^2 \leq 0\,\,,
$$
so $M$ is a decreasing quantity. Also, since the profile $\psi$ corresponds to a smooth solution, one has necessarily 
$\psi(0) = 0$ (modulo the symmetries of the equation), so $M(0)=0$. 

Now, arguing by contradiction, suppose that $\psi(1) = 0$, then $M(1) = 0$ and $M$ has to be identically $0$ on $[0,1]$. By the formula for $M'$, this implies that $\psi$ is constant on $[0,1]$, which is a contradiction.

\bigskip

So $\psi$ satisfies $\psi (1) = \phi^*$, and if we define
$$
\bar{\psi}(\rho) = \left\{ \begin{array}{ll} \psi(\rho) & \mbox{if $\rho \leq 1$} \\
\phi^* & \mbox{if $\rho \geq 1$}\,\,,
\end{array} \right.
$$
the associated space-time function
$$
\phi(r,t) = \psi\left( \frac{r}{t} \right)
$$
solves $(EWM)$ with data equal to the equator map, though it is not equal to the equator map itself. Thus
$(ii)$ does not hold true, and neither does $(i)$ by the previous section.

\section{Proof of Proposition~\ref{energyconserve}}

\label{pec}

We are considering a smooth profile $\psi$, such that $\psi(0) =0$ and $\psi(1)=\phi^*$, and
$$
\phi(r,t) = \left\{ \begin{array}{ll} \psi \left( \frac{r}{t} \right) & \mbox{if $r\leq t$ }\\  \phi^* & \mbox{if $r \geq t$\,\,.} \end{array} \right.
$$
solves $(EWM)$. Our aim is to show that $\phi$ satisfies the energy equality~(\ref{energyeq}). 

We begin with the

\begin{claim}
There holds $\displaystyle \frac{g(\phi)g'(\phi)}{r^2} \in L^1_{t\;\;loc} L^2(r^{d-1}dr)$.
\end{claim}

\textsc{Proof:} Our function is non zero only for $r<t$. Thus we compute
\begin{equation*}
\begin{split}
\int_0^t \left[ \frac{g(\phi)g'(\phi)}{r^2} \right]^2 r^{d-1} dr & \leq C \int_0^t \left| \frac{\psi \left( \frac{r}{t} \right) }{r^2} \right|^2 r^{d-1} dr \\
& \leq C t^{d-4} \int_0^1 \left| \frac{\psi(s)}{s^2} \right|^2 s^{d-1} ds \\
& \leq C t^{d-4} \int_0^1 s^{d-3} \,ds \\
& \leq C t^{d-4} \,\,.
\end{split}
\end{equation*}
Thus
$$
\left\| \frac{g(\phi)g'(\phi)}{r^2} \right\|_{L^2(r^{d-1}dr)} \leq C t^{\frac{d-4}{2}} \,\,,
$$
which is a locally integrable function of $t$ for $d \geq 3$. $\blacksquare$

\bigskip

The proof of the energy equality is now standard. Define the mollification operator by the $\mathcal{C}_0^\infty$ function $Z$,
$$
f^\epsilon (x) = \left[ f * \frac{1}{\epsilon^d} Z \left( \frac{\cdot}{\epsilon} \right) \right] (x)\,\,.
$$
Notice that this operator is defined for general functions of $x$, but we can apply it to radial functions, thus depending only on $r$, too.

Now we mollify the equation, and take the $L^2(r^{d-1}dr)$ scalar product with $\phi^\epsilon_t$ over the truncated cone given by
$$
K_{0,T,R} = \{ (r,t) \;\;\mbox{such that $0 \leq t \leq T$ and $0 \leq r \leq R - t$} \} \,\,.
$$
Thus we have
\begin{equation}
\label{orchidee}
\int_K \left[ \phi_{tt} -\phi_{rr} - \frac{d-1}{r} \phi_r + \frac{k}{r^2} g(\phi) g'(\phi) \right]^\epsilon \phi^\epsilon_t \,r^{d-1}dr\,dt = 0 \,\,.
\end{equation}
Of course the idea is to let $\epsilon$ go to zero and recover the energy equality~(\ref{energyeq}). 

It is easily seen that
\begin{equation*}
\begin{split}
\int_K & \left[ \phi_{tt} -\phi_{rr} - \frac{d-1}{r} \phi_r \right]^\epsilon \phi^\epsilon_t r^{d-1}dr\,dt \overset{\epsilon \rightarrow 0}{\longrightarrow} \int_0^{R-T} 
\left( \phi_r^2(T) + \phi_t^2(T) \right) r^{d-1} \,dr \\
&\;\;\;\;\;\;\;\;\;\;\;\;\;\;\;\;\;\;\;\;\;\;\;\;\;\;\;\;\;\;\;\;\;\;\;\;\;\;\;\;\;\;\; - \int_0^R \left( \phi_r^2(0) + \phi_t^2(0) \right) r^{d-1} \,dr \\
&\;\;\;\;\;\;\;\;\;\;\;\;\;\;\;\;\;\;\;\;\;\;\;\;\;\;\;\;\;\;\;\;\;\;\;\;\;\;\;\;\;\;\; + \int_{S_{0,R,T}} \left( \frac{1}{2} \phi_t^2 + \frac{1}{2} |\nabla \phi|^2 - \phi_t \nabla \phi \cdot n \right) d\sigma \,\,.
\end{split}
\end{equation*}
Only the last summand in~(\ref{orchidee}) is problematic, and we rewrite it as
$$
\int_{K \cap \{ r \geq \delta\} } \left[ \frac{k}{r^2} g(\phi) g'(\phi) \right]^\epsilon \phi^\epsilon_t r^{d-1}dr\,dt + \int_{K \cap \{ r \leq \delta \}} \left[ \frac{k}{r^2} g(\phi) g'(\phi) \right]^\epsilon \phi^\epsilon_t r^{d-1}dr\,dt
\overset{def}{=} I_{\delta,\epsilon} + II_{\delta, \epsilon} \,\,.
$$
The conclusion follows since
\begin{equation*}
\begin{split}
I_{\delta,\epsilon} & \overset{\epsilon \rightarrow 0}{\longrightarrow} 
\int_\delta^{R-T} \frac{k}{r^2} g(\phi(T))^2 r^{d-1} \,dr
- \int_\delta^{R} \frac{k}{r^2} g(\phi(0))^2 r^{d-1} \,dr \\
&\;\;\;\;\;\;\;\;\;\;\;\;\;\;\;\;\;\;\;\;\;\;\;\;\;\;\;\;\;\;\;\;\;+ \int_{S_{0,R,T} \setminus \{r \leq \delta \}} \frac{k}{2r^2} g(\phi)^2 d\sigma \\
& \overset{\delta \rightarrow 0}{\longrightarrow} \int_0^{R-T} \frac{k}{r^2} g(\phi(T))^2 r^{d-1} \,dr
- \int_0^{R} \frac{k}{r^2} g(\phi(0))^2 r^{d-1} \,dr
+ \int_{S_{0,R,T}} \frac{k}{2r^2} g(\phi)^2 d\sigma \\
\end{split}
\end{equation*}
and
$$
II_{\delta,\epsilon} \overset{\epsilon \rightarrow 0}{\longrightarrow} \int_{K \cap \{ r \leq \delta \}}  \frac{k}{r^2} g(\phi) g'(\phi) \phi_t r^{d-1}dr\,dt \overset{\delta \rightarrow 0}{\longrightarrow} 0 \,\,.
$$

\section{Regularity of the profile}

\label{proofreg}

\subsection{Profiles in $\dot{H}^{3/2}$ }

In this subsection, we prove that if $R_\chi \psi \in \dot{H}^{3/2}(B_{\mathbb{R}^3}(0,1))$ and $\psi$ solves~(\ref{eqpsi}), then $R_\chi \psi \in \mathcal{C}^\infty (\bar B_{\mathbb{R}^d}(0,1))$.

It is clear that $\psi \in \mathcal{C}^\infty ((0,1))$ so we simply need to examine the points $0$ and $1$.

\subsubsection{Smoothness of $\psi$ around $0$}

Define
$$
M(\rho) = \frac{1}{2} (\rho^2 - \rho^4) \psi'(\rho)^2 - k g(\psi(\rho))^2
$$
and observe that
$$
M'(\rho) \leq 0\,\,,
$$
hence $M$ is decreasing.

\bigskip

Assume first that $M$ is bounded at $0$.

Then, as is proved in~\cite{Ge}, $\psi$ is continuous at $0$ and $\psi(0) = 0$ or $\phi^*$. 

Since $M$ is decreasing, and $g$ is maximal at $\phi^*$, $\psi(0) = \phi^*$ is possible only if $\phi$ is identically equal to $\phi^*$. But then $R_\chi \psi$ does not belong to $\dot{H}^{3/2}$.

Thus necessarily $\psi(0) = 0$. Then $R_\chi \psi$ is continuous at $0$, and it satisfies the harmonic map equation from the hyperbolic ball to $N$. By elliptic regularity (see for instance~\cite{LU}), $R_\chi \psi$ is smooth in a neighborhood of $0$.

\bigskip

We now assume that $M$ is not bounded at $0$, and we will reach a contradiction. Since $g$ is bounded, and $M$ decreasing, $M$ can be unbounded at $0$ if and only if for a constant $C$ and close to $0$
$$
|\psi'(\rho)| \geq \frac{C}{\rho} \,\,.
$$
But this implies that
$$
\partial_r R_\chi \psi = ( \psi'(r) \,,\,0) \notin L^3 (B_{\mathbb{R}^3}(0,1))
$$
(the coordinates correspond to the tangent space to $N$ at $\psi$). This is a contradiction since 
\begin{equation*}
\begin{split}
R_\chi \psi \in \dot{H}^{3/2} (B_{\mathbb{R}^3}(0,1)) & 
\;\;\;\Rightarrow\;\;\; \nabla_x R_\chi \psi \in \dot{H}^{1/2} (B_{\mathbb{R}^3}(0,1)) \hookrightarrow L^3(B_{\mathbb{R}^3}(0,1)) \\
& \Rightarrow\;\;\; \partial_r R_\chi \psi = \frac{x}{|x|} \cdot \nabla_x R_\chi \psi \in L^3(B_{\mathbb{R}^3}(0,1))\,\,.
\end{split}
\end{equation*}

\subsubsection{Smoothness of $\psi$ near $1$}

We notice first that if $R_\chi \psi \in \dot{H}^{3/2}$, $\psi$ is continuous at $1$ thus $\psi(1)$ is well-defined.

If $\psi(1) = \phi^*$, regularity can be proved as in~\cite{STZ}, Lemma 4.2.

If $\psi(1) = 0$, one can show as in Subsection~\ref{subsubiii} that $\psi$ is constant.

Thus the only remaining case is
$$
g(\psi(1)) g'(\psi(1)) \neq 0 \,\,.
$$
We will assume that this holds and reach a contradiction.

Switching to the unknown function
$$
\widetilde{\psi} (z) = \psi \left(\frac{1}{z} \right) \,\,,
$$
we get the new equation
\begin{equation}
\label{eqpsitilde}
(z^2 - 1) \widetilde{\psi}'' (z) = g(\widetilde \psi(z)) g'(\widetilde \psi(z))
\end{equation}
Integrating it twice yields successively
\begin{equation*}
\begin{split}
& \widetilde{\psi}' (z) \sim g(\psi(1)) g'(\psi(1)) \log|z-1| \;\;\;\;\mbox{as $z \rightarrow 1$} \\
& \widetilde{\psi}(z) = \psi(1) + O(|z-1|\log|z-1|) \;\;\;\;\mbox{as $z \rightarrow 1$} \,\,.
\end{split}
\end{equation*}
Inserting this development in~(\ref{eqpsitilde}) gives
$$
\widetilde{\psi}'' (z) = \frac{g(\psi(1)) g'(\psi(1))}{z^2 - 1} + O (\log |z-1|) \,\,.
$$
Integrating this new equation yields
$$
\widetilde{\psi}' (z) = C + \frac{ g(\psi(1)) g'(\psi(1))}{2} \log |z-1| + D(z)
$$
where $C$ is a constant and $|D'(z)| \leq |\log|z-1||$.

Since we are considering an equivariant problem around $r =1$, the regularity can be studied by reducing matters to the 1-dimensional case. It suffices to observe that
$$
( z \mapsto D(z) ) \in \dot{H}^{1/2}_{loc}(\mathbb{R})\,\,,
$$
but on the other hand
$$
z \mapsto \log|z| \in \dot{B}^{1/2}_{2,\infty \;\; loc}(\mathbb{R}) \setminus \dot{H}^{1/2}_{loc}(\mathbb{R})
$$
(this is a simple computation, done in~\cite{Ge}) to conclude that
\begin{equation}
\label{castor}
\widetilde{\psi}' (z) \in \dot{B}^{1/2}_{2,\infty \;\; loc}(\mathbb{R}) \setminus \dot{H}^{1/2}_{loc}(\mathbb{R})
\end{equation}
which is the desired contradiction.

\subsection{Profiles in $\dot{B}^{3/2}_{2,\infty}$}

Our claim is that without any further assumptions than the ones of Section~\ref{lezard}, there exist non-constant profiles such that
$$
R_\chi \psi \in \dot{B}^{3/2}_{2,\infty}(\mathbb{R}^3)
$$
and $R_\chi \psi$ solves $(WM)$.

\bigskip

We sketch the proof, which simply consists of solving the ODE
$$
\left\{ \begin{array}{l} \psi'' + \frac{2}{\rho} \psi' + \frac{k g(\psi) g'(\psi)}{\rho^2(1-\rho^2)} = 0 \\
\psi(0) = \dots = \psi^{(l-1)}(0) = 0 \\ \psi^{(l)}(0) = \alpha \end{array} \right.
$$
(see~\cite{CSTZ}, Lemma 2.5). One gets a solution which is smooth on $[0,1)$. We see as in~(\ref{castor}) that $ \psi $ is in $\dot{B}^{3/2}_{2,\infty}$ around $1$. This concludes the proof.

\section{Proof of Proposition~\ref{PGProp}}

\label{sectionODE}

Let us recall the statement of Proposition~\ref{PGProp}: we aim at proving that if $\psi$ is locally integrable, and smooth at $0$,
$$
u(x,t) = R_\chi \psi \left( \frac{|x|}{t} \right)
$$
solves the wave-map equation if and only if $\psi$ solves
\begin{equation}
\label{ODEpsi}
\rho^2 (1-\rho^2) \psi''(\rho) + \rho ((d-1)- 2 \rho^2) \psi'(\rho) - k g(\psi) g'(\psi) = 0
\end{equation}
classically on $(0,1) \cup (1,\infty)$ and furthermore
$$
\psi(1^-) = \psi(1^+)\,\,.
$$

\subsection{From the PDE to the ODE}

A first step consists in proving that $u$ as defined above solves $(WM)$ for $t>0$ if and only if $\psi$ solves~(\ref{ODEpsi}) in $\mathcal{S}'$. This can be done in a straightforward manner, we refer to~\cite{Ge} for a proof in a less general situation.

What remains is to see that the above ODE is solved in $\mathcal{S}'$ if and only if it is solved classically on $(0,1) \cup (1,\infty)$ and furthermore $\psi(1^-) = \psi(1^+)$. This is the aim of the next subsection.

\subsection{Condition for solving the ODE in $\mathcal{S}'(0,1)$}

Solving~(\ref{ODEpsi}) in $\mathcal{S}'$ or in the classical sense is equivalent in $(0,1)$ and $(1,\infty)$. We now examine the situation of $\rho =1$.

\bigskip

The equation can be written away from 1 as
$$
\frac{d}{d\rho} \left( \frac{\rho^{d-1}}{(\rho^2-1)^{\frac{d-3}{2}}} \psi' \right) = - \frac{\rho^{d-3}}{(\rho^2 -1)^{\frac{d-1}{2}}} k g(\psi) g'(\psi) \,\,.
$$
Integrating it starting from say $1/2$ or $3/2$, to $\rho$, gives
\begin{equation}
\label{estpsi}
\left| \frac{1}{(\rho-1)^{\frac{d-3}{2}}} \psi'(\rho) \right| \leq C +
\left\{
\begin{array}{l} C |\log(\rho-1)| \;\;\;\;\mbox{if $d=3$} \\
\frac{C}{|\rho-1|^{\frac{d-3}{2}}} \;\;\;\;\mbox{if $d \geq 4$} \,\,.
\end{array}
\right.
\end{equation}
In any case, $\psi'$ is integrable on $(0,1) \cup (1,\infty)$, and $\psi$ is left and right-continuous at $1$. Set
$$
a \overset{def}{=} \psi(1^+) - \psi(1^-)\,\,.
$$
$\psi'$ can be written as
\begin{equation}
\label{laitue}
\psi' = \xi + a \delta_1 \,\,,
\end{equation}
where $\xi$ belongs to $L^1$ and $\delta_1$ is a Dirac weight at $\rho=1$.

\bigskip

Since~(\ref{ODEpsi}) holds in the classical sense away from $1$, we see that it holds in $\mathcal{S}'(0,\infty)$ if and only if for any $f \in \mathcal{S}$,
$$
\lim_{\epsilon \rightarrow 0} \int_{1-\epsilon}^{1+\epsilon} \left[ \rho^2 (1-\rho^2) \psi''(\rho) + \rho ((d-1)- 2 \rho^2) \psi'(\rho) - k g(\psi) g'(\psi) \right] f(\rho) \,d\rho = 0
$$
(where the integral is to be understood in the distribution sense). Using the previous discussion, this can be rewritten as
$$
(d-3)a f(1) + \lim_{\epsilon \rightarrow 0} \int_{1-\epsilon}^{1+\epsilon} \rho^2 (1-\rho^2) \psi''(\rho) f(\rho) \,d\rho = 0
$$
or, integrating by parts,
$$
(d-3)a f(1) +
\lim_{\epsilon \rightarrow 0}  \left[ \psi' \rho^2 (1-\rho^2) f \right]_{1-\epsilon}^{1+\epsilon} - \int_{1-\epsilon}^{1+\epsilon} \psi' \frac{d}{d \rho} ( \rho^2 (1-\rho^2) f ) \,d\rho = 0 \,\,
$$
which, by the estimates~(\ref{estpsi}) and~(\cite{laitue}), is equivalent to
$$
(d-1)a f(1) = 0 \,\,.
$$
Thus we get $a =0$, which is the desired result.

\bigskip

{\bf Acknowledgements:} The author is indebted to Jalal Shatah for very interesting discussions during the writing of this article.

\end{document}